# CONDITIONING AND INITIAL ENLARGEMENT OF FILTRATION ON A RIEMANNIAN MANIFOLD[1]

By Fabrice Baudoin

*Université Paul Sabatier*

We extend to Riemannian manifolds the theory of conditioned stochastic differential equations. We also provide some enlargement formulas for the Brownian filtration in this nonflat setting.

**1. Introduction.** In this paper we develop the theory of conditioned stochastic differential equations (CSDEs) on a Riemannian manifold. In the flat case, this theory was initiated in [1] and further used in [3] and [4]. In [1] and [3], we gave an application of the CSDEs to the mathematical finance topic of informed insiders. In [4], from the drift of the CSDEs we first constructed some martingales (called Newton's martingales) which generalize the stochastic Newton equation for the so-called reciprocal processes (see [28]). Then we studied the symmetries of the CSDEs, that is, the transformations on the flat path space which preserve the set of CSDEs constructed from a given functional. From this we constructed some martingales called Noether's martingales, by analogy with the classical Noether theorem. This paper can be read independently of the articles cited above.

In the present paper, we study a Riemannian Brownian motion with drift $V$ for which one functional of the trajectories is conditioned (in Doob's sense) to follow a given law. The conditioned process, which generalizes naturally the conditioned diffusion (see [5]) that Bismut used extensively in his probabilistic proof of Atiyah–Singer theorems (see [6, 7]), is shown to be a semimartingale in its own filtration. It is also shown that it is semimartingale in the Brownian filtration initially enlarged by the conditioned functional until each time smaller than the revelation time of the functional. Furthermore, this decomposition does not depend on the law of the conditioning. Then we study the case of conditioning of a marginal law for Riemannian Brownian

---

Received October 2002; revised August 2003.

[1]Supported by the Research Training Network HPRN-CT-2002-00281.

*AMS 2000 subject classifications.* 58G32, 60J60, 60H10.

*Key words and phrases.* Conditioned stochastic differential equations, diffusions on manifolds, initial enlargement of filtration, Malliavin–Bismut calculus of variations.

motion with drift. Finally, we are interested in the conditioning of hitting times for small geodesic spheres, and we give a probabilistic characterization of rank 1 symmetric spaces which is derived from [2] and [26].

As will be seen, actually, it appears that most of the results presented in [1] can be extended to Riemannian manifolds. Roughly speaking, this extension can be explained by the fact that the horizontal lifting commutes with the conditioning. In other words, if we condition a Brownian motion by a functional and if we lift this conditioned process, then we obtain a horizontal Brownian motion conditioned by the lifted functional.

**2. Framework and assumptions.** We now turn to the notations which are used throughout the paper. Let $(\mathbb{M}, g)$ be a $d$-dimensional connected complete Riemannian manifold. We denote by $\Delta$ the Laplace–Beltrami operator on $\mathbb{M}$ (for us, $\Delta$ is negative). The tangent bundle to $\mathbb{M}$ is denoted $T\mathbb{M}$ and $T_m\mathbb{M}$ is the tangent space at $m$: we have, hence, $T\mathbb{M} = \bigcup_n T_n\mathbb{M}$. The orthonormal frame bundle of $\mathbb{M}$ is denoted $\mathcal{O}(\mathbb{M})$. Hence, $(\mathcal{O}(\mathbb{M}), \mathbb{M}, \mathcal{O}_d(\mathbb{R}))$ is a principal bundle on $\mathbb{M}$ with structure group $\mathcal{O}_d(\mathbb{R})$ of $d \times d$ orthogonal matrices. The transpose of a matrix $M$ is denoted $^TM$. We denote by $\pi$ the canonical surjection $\mathcal{O}(\mathbb{M}) \to \mathbb{M}$. The horizontal fundamental vector fields of $\mathcal{O}(\mathbb{M})$ are denoted $(H_i)_{i=1,\dots,d}$. The Bochner horizontal Laplacian, that is, the lift of $\Delta$, is then given by

$$\Delta_{\mathcal{O}(\mathbb{M})} = \sum_{i=1}^{d} H_i^2.$$

The symbol $\nabla$ denotes the covariant differentiation on $\mathbb{M}$ associated with the torsion-free connection on $\mathbb{M}$ (the Levi–Civita connection). In addition, **Ric** denotes the Ricci curvature tensor and $\overline{\mathbf{Ric}}$ denotes its equivariant representation. For $u \in \mathcal{O}(\mathbb{M})$, $\overline{\mathbf{Ric}}_u$ is hence an application $\mathbb{R}^d \to \mathbb{R}^d$. For a smooth vector field $V$ on $\mathbb{M}$, the equivariant representation of the $(1,1)$ tensor $\nabla V$ is denoted $\overline{\nabla V}$, and $\overline{\nabla V}(u)$ is thus also an application $\mathbb{R}^d \to \mathbb{R}^d$.

Let us now consider $m \in \mathbb{M}$, and a smooth vector field $V$. We associate the stochastic differential equation on $\mathcal{O}(\mathbb{M})$,

$$Z_t^* = U_0 + \int_0^t V^*(Z_s^*)\, ds + \sum_{i=1}^{d} \int_0^t H_i(Z_s^*) \circ dB_s^i,$$

where:

- $U_0 \in \mathcal{O}(\mathbb{M})$ is such that $\pi U_0 = m$;
- $V^*$ is the lift of $V$;
- $\circ$ denotes the integration in a Stratonovitch sense;
- $(B_t)_{t \geq 0}$ is a $d$-dimensional standard Brownian motion.

We assume that this equation has, for any $m \in \mathbb{M}$, a unique strong non-explosive solution in the sense that there exists on $\mathcal{O}(\mathbb{M})$ a unique process $(Z_t^*)_{t \geq 0}$ whose natural filtration is equal to the natural filtration of $(B_t)_{t \geq 0}$ and such that for any $C^\infty$ bounded function $f : \mathcal{O}(\mathbb{M}) \to \mathbb{R}$,

$$f(Z_t^*) = f(U_0) + \int_0^t (V^*f)(Z_s^*)\, ds + \sum_{i=1}^d \int_0^t (H_i f)(Z_s^*) \circ dB_s^i, \qquad t \geq 0.$$

Consider the space of continuous paths

$$\mathcal{C}_m(\mathbb{M}) = \{\omega : \mathbb{R}_+ \to \mathbb{M}, \ \omega(0) = m, \ \omega \text{ continuous}\}.$$

From the previous assumption there exists a unique probability measure $\mathbb{P}_m$ on $\mathcal{C}_m(\mathbb{M})$ such that for any $C^\infty$ bounded $f : \mathbb{M} \to \mathbb{R}$, the process

$$(2.1) \qquad \left( f(X_t) - f(m) - \frac{1}{2} \int_0^t (\Delta f)(X_s)\, ds - \int_0^t (Vf)(X_s)\, ds \right)_{t \geq 0}$$

is, under $\mathbb{P}_m$, an $\mathcal{F}$-adapted martingale null at 0. Here $(X_t)_{t \geq 0}$ denotes the coordinate process on $\mathcal{C}_m(\mathbb{M})$ and $(\mathcal{F}_t)_{t \geq 0}$ denotes its natural filtration. In other words, the law of the process $Z = \pi Z^*$ is the unique solution of the martingale problem with initial condition $m$ associated with the elliptic operator

$$\tfrac{1}{2}\Delta + V.$$

We refer to [24] and [25] for the general theory of the so-called martingale problems.

The transition function of $(Z_t)_{t \geq 0}$ is denoted $q_t$, hence we have, for $s < t$,

$$\mathbb{P}(Z_t \in dy | \mathcal{F}_s) = q_{t-s}(Z_s, y)\, dy,$$

where $dy$ is the Riemannian volume measure on $\mathbb{M}$. The existence of the function $q_t$ and its smoothness comes from Hörmander's theorem. Moreover, we assume that $q_t$ is positive.

**3. Conditioning and pinning class in the nonflat Wiener space.** We fix now once and for all $m$ and $U_0$. For the sake of simplicity, we denote simply by $\mathbb{P}$ the probability measure $\mathbb{P}_m$ on $\mathcal{C}_m(\mathbb{M})$ described in the previous section. Since we mainly focus our attention on the laws of $\mathbb{M}$-valued processes, we work in the stochastic basis

$$(\mathcal{C}_m(\mathbb{M}), (X_t)_{t \geq 0}, (\mathcal{F}_t)_{t \geq 0}, \mathbb{P}).$$

For this, we have to transfer the assumptions of the previous section into this stochastic basis. Namely, there are a unique process $(X_t^*)_{t \geq 0}$ on $\mathcal{O}(\mathbb{M})$ and a unique $d$-dimensional standard Brownian motion $(\tilde{X}_t)_{t \geq 0}$ such that

$$X_t^* = U_0 + \int_0^t V^*(X_s^*)\, ds + \sum_{i=1}^d \int_0^t H_i(X_s^*) \circ d\tilde{X}_s^i, \qquad t \geq 0, \ \pi X^* = X.$$

We consider now on $\mathcal{C}_m(\mathbb{M})$ a random variable $Y$-valued in some Polish space $\mathcal{P}$, endowed with Borel $\sigma$-algebra $\mathcal{B}(\mathcal{P})$, and measurable with respect to the $\sigma$-algebra $\mathcal{F}_T$ with $T \in \mathbb{R}_+ \cup \{+\infty\}$. We assume the existence of a regular disintegration of $Y$ with respect to the filtration $\mathcal{F}$. Namely, we assume that there exists a jointly measurable, continuous in $t$ and $\mathcal{F}$-adapted process

$$(3.1) \qquad \eta_t^y, \qquad 0 \leq t < T, \; y \in \mathcal{P},$$

satisfying, for $dt \otimes \mathbb{P}_Y$ almost every $0 \leq t < T$ and $y \in \mathcal{P}$,

$$\mathbb{P}(Y \in dy | \mathcal{F}_t) = \eta_t^y \mathbb{P}_Y(dy),$$

where $\mathbb{P}_Y$ denotes the law of $Y$ under $\mathbb{P}$.

REMARK 1. For the sake of presentation, we restrict ourselves to the case of a deterministic horizon time $T$. Nevertheless the contents contained in this section and the next one (excepted what deals with Malliavin calculus) are easily extended by taking for $T$ a stopping time of the filtration $\mathcal{F}$ and taking for $Y$ an $\mathcal{F}_T$-measurable functional. In that case, we have to work under the assumption that there exists a jointly measurable process $(\eta_t^y, 0 \leq t < T, y \in \mathcal{P})$ satisfying, for all bounded and measurable function $f$,

$$\mathbb{E}(f(Y) | \mathcal{F}_t, \, t < T) = \int_{\mathcal{P}} f(y) \eta_t^y \mathbb{P}_Y(dy).$$

Actually, to include the case of a random horizon in our presentation, it would suffice to work in the filtration $(\mathcal{F}_t \cap \{t < T\})_{t \geq 0}$ (such a "trick" is well known in the theory of enlargement of filtrations and can, e.g., be found in [29]).

One of our main objects of study is the so-called pinning class (see [4]) of the measure $\mathbb{P}$ with respect to the functional $Y$, that is, the set $\mathcal{R}_Y(\mathbb{P})$ of probability measures on $\mathcal{C}_m(\mathbb{M})$ defined by

$$\mathcal{R}_Y(\mathbb{P}) = \{\mathbb{Q} \sim \mathbb{P}, \mathbb{Q}(\cdot | Y) = \mathbb{P}(\cdot | Y)\}.$$

To explicate the semimartingale decomposition of $X$ under $\mathbb{Q} \in \mathcal{R}_Y(\mathbb{P})$, we need the following nonflat version of the so-called Jacod lemma in the theory of initial enlargement of filtration (see [18]). In what follows, $\mathcal{P}(\mathcal{F})$ denotes the predictable $\sigma$-field associated with the filtration $\mathcal{F}$.

LEMMA 2. There exists a $\mathcal{P}(\mathcal{F}) \otimes \mathcal{B}(\mathcal{P})$ measurable process

$$[0, T[ \times \mathcal{C}_m(\mathbb{M}) \times \mathcal{P} \to T\mathbb{M},$$

$$(t, \omega, y) \to \alpha_t^y(\omega)$$

such that:

1. *For $\mathbb{P}_Y$-a.e. $y \in \mathcal{P}$ and for $0 \leq t < T$, $\alpha_t^y \in T_{X_t}\mathbb{M}$.*
2. *For $\mathbb{P}_Y$-a.e. $y \in \mathcal{P}$ and for $0 \leq t < T$,*

$$\mathbb{P}\left(\int_0^t \|\alpha_u^y\|^2 \, du < +\infty\right) = 1.$$

3. *For $\mathbb{P}_Y$-a.e. $y \in \mathcal{P}$ and for all 1-forms $\theta$ on $\mathbb{M}$,*

$$\left\langle \eta^y, \int_{X[0,\cdot]} \theta \right\rangle_t = \int_0^t \eta_u^y \theta(\alpha_u^y) \, du, \qquad 0 \leq t < T.$$

Proof.    From the predictable representation property of the Brownian motion $(\tilde{X}_t)_{t \geq 0}$ and from [18], there exists a $\mathcal{P}(\mathcal{F}) \otimes \mathcal{B}(\mathcal{P})$ measurable process

$$[0, T[ \times \mathcal{C}_m(\mathbb{M}) \times \mathcal{P} \to \mathbb{R}^d,$$
$$(t, \omega, y) \to \tilde{\alpha}_t^y(\omega)$$

such that:

1. For $\mathbb{P}_Y$-a.e. $y \in \mathcal{P}$ and for $0 \leq t < T$,

$$\mathbb{P}\left(\int_0^t \|\tilde{\alpha}_u^y\|^2 \, du < +\infty\right) = 1.$$

2. For $\mathbb{P}_Y$-a.e. $y \in \mathcal{P}$ and for $0 \leq t < T$, $1 \leq i \leq d$,

$$\langle \eta^y, \tilde{X}^i \rangle_t = \int_0^t \eta_u^y \tilde{\alpha}_u^{y,i} \, du.$$

Then we set

$$\alpha_t^y = X_t^* \tilde{\alpha}_t^y$$

and it is easy to verify that it satisfies the conditions of the lemma. In particular, let us show that it satisfies the fourth condition. For a 1-form $\theta$ on $\mathbb{M}$, we have

$$\left\langle \eta^y, \int_{X[0,\cdot]} \theta \right\rangle_t = \left\langle \eta^y, \sum_{i=1}^d \int_0^\cdot \theta(X_s^* e_i) \circ d\tilde{X}_s^i \right\rangle_t,$$

where $(e_i)_{i=1,\dots,d}$ is the canonical base of $\mathbb{R}^d$. However,

$$\left\langle \eta^y, \int_0^\cdot \theta(X_s^* e_i) \circ d\tilde{X}_s^i \right\rangle_t = \int_0^t \eta_s^y \theta(X_s^* e_i) \tilde{\alpha}_s^{y,i} \, ds.$$

Thus,

$$\left\langle \eta^y, \int_{X[0,\cdot]} \theta \right\rangle_t = \int_0^t \eta_u^y \theta(\alpha_u^y) \, du, \qquad 0 \leq t < T. \qquad \square$$

On $\alpha$ we furthermore make the following integrability assumption: For almost every $t \in [0, T)$, it holds $\mathbb{P}$-almost surely that

$$(3.2) \qquad \int_0^t \mathbb{E}(\|\alpha_u^Y\| | \mathcal{F}_u)^2 \, du < +\infty.$$

We can now deduce the following proposition.

PROPOSITION 3. *Let $\mathbb{Q} \in \mathcal{R}_Y(\mathbb{P})$. Then under $\mathbb{Q}$, the coordinate process $(X_t)_{0 \leq t \leq T}$ is a semimartingale in the filtration $(\mathcal{F}_t)_{0 \leq t \leq T}$. Moreover, the process*

$$\tilde{X}_t - \int_0^t (X_s^*)^{-1} \frac{\int_\mathcal{P} \eta_s^y \alpha_s^y \mathbb{Q}_Y(dy)}{\int_\mathcal{P} \eta_s^y \mathbb{Q}_Y(dy)} \, ds, \qquad t < T,$$

*is a Brownian motion under $\mathbb{Q}$, where $\mathbb{Q}_Y$ is the law of $Y$ under $\mathbb{Q}$.*

PROOF. Let $\mathbb{Q} \in \mathcal{R}_Y(\mathbb{P})$ and denote by $\xi$ the Radon–Nikodym density $d\mathbb{Q}_Y/d\mathbb{P}_Y$ which is well defined because $\mathbb{Q} \sim \mathbb{P}$. We have, for any $\mathcal{F}_T$-measurable positive and bounded random variable $F$,

$$\mathbb{E}^\mathbb{Q}(F) = \int_\mathcal{P} \mathbb{E}^\mathbb{Q}(F | Y = y) \xi(y) \mathbb{P}_Y(dy).$$

However, since $\mathbb{Q} \in \mathcal{R}_Y(\mathbb{P})$ and

$$\int_\mathcal{P} \mathbb{E}^\mathbb{Q}(F | Y = y) \xi(y) \mathbb{P}_Y(dy) = \int_\mathcal{P} \mathbb{E}(F | Y = y) \xi(y) \mathbb{P}_Y(dy) = \mathbb{E}(\xi(Y) F),$$

we conclude that

$$\mathbb{Q}_{/\mathcal{F}_T} = \xi(Y) \mathbb{P}_{/\mathcal{F}_T}.$$

From this and from

$$\mathbb{P}(Y \in dy | \mathcal{F}_t) = \eta_t^y \mathbb{P}_Y(dy)$$

we deduce that for $t < T$,

$$\mathbb{Q}_{/\mathcal{F}_t} = \left( \int_\mathcal{P} \eta_t^y \xi(y) \mathbb{P}_Y(dy) \right) \mathbb{P}_{\mathcal{F}_t} = \left( \int_\mathcal{P} \eta_t^y \mathbb{Q}_Y(dy) \right) \mathbb{P}_{\mathcal{F}_t}.$$

Now, we can conclude with Girsanov's theorem, because from the proof of Lemma 1, for $\mathbb{P}_Y$-a.e. $y \in \mathcal{P}$ and for $0 \leq t < T$, $1 \leq i \leq d$,

$$\langle \eta^y, \tilde{X}^i \rangle_t = \int_0^t \eta_u^y \tilde{\alpha}_u^{y,i} \, du$$

with

$$\tilde{\alpha}_t^y = (X_t^*)^{-1} \alpha_t^y.$$

Because of the assumption (3.2) we can apply Fubini's theorem

$$\left\langle \int_{\mathcal{P}} \eta^y_{\cdot} \mathbb{Q}_Y(dy), \tilde{X}^i \right\rangle_t = \int_0^t \int_{\mathcal{P}} \eta^y_u \tilde{\alpha}^{y,i}_u \mathbb{Q}_Y(dy)\, du,$$

which leads to the expected result. $\square$

In the next theorem, under further regularity assumptions we try to compute more explicitly the compensator of $\tilde{X}$ under $\mathbb{Q} \in \mathcal{R}_Y(\mathbb{P})$ by means of the Clark–Ocone–Bismut formula (see [15]; [21], Theorem 6.4, and [5], Theorem 2.2, pages 61 and 62). In particular, we see a Bakry–Emery curvature type term which appears in the computations and which measures exactly the difference with the flat case. This term stems from Weitzenböck formula on 1-forms. Before we state our formula, let us recall some basic facts about Malliavin calculus on a nonflat space. For a cylindrical functional,

$$F = f(X_{t_1}, \ldots, X_{t_n}),$$

where $f \colon \mathbb{M}^n \to \mathbb{R}$ is a smooth function and the directional derivative of $F$ along the Cameron–Martin vector field $D_h$ is given by

$$D_h F = \sum_{i=1}^n \left(\nabla^i f(X_{t_1}, \ldots, X_{t_n}), X^*_{t_i} h\right),$$

where $h$ is an $\mathbb{R}^d$-valued adapted process with a derivative in $L^2$ such that

$$\mathbb{E}\left(\int_0^T \left(\frac{dh}{dt}\right)^2\right) < +\infty.$$

Now, the Malliavin derivative of $F$ is defined by the representation formula

$$D_h F = \int_0^T \left(\mathbf{D}_s F, \frac{dh}{ds}\right) ds.$$

It is shown in the same way as in the flat case (see [22], page 26) that $D_h \colon \mathcal{S} \to L^p(\mathcal{C}_m(\mathbb{M}), \mathbb{P})$ is closable for $p \geq 1$, where $\mathcal{S}$ is the set of cylindrical functionals.

Finally, we denote by $\mathbf{\Omega}$ the field of linear applications $\mathbb{R}^d \to \mathbb{R}^d$ defined by

$$\mathbf{\Omega} := \tfrac{1}{2} \overline{\mathbf{Ric}} -^T \overline{\nabla V}$$

and we assume that it is bounded.

PROPOSITION 4. *Let $\mathbb{Q} \in \mathcal{R}_Y(\mathbb{P})$ and assume, moreover, that $\xi := d\mathbb{Q}_Y/d\mathbb{P}_Y$ has a version such that $\xi(Y)$ and $\ln\xi(Y) \in \mathrm{Dom}(\mathbf{D})$. Then*

$$(3.3)\quad \begin{aligned} &(X^*_t)^{-1}\frac{\int_{\mathcal{P}} \eta^y_t \alpha^y_t \mathbb{Q}_Y(dy)}{\int_{\mathcal{P}} \eta^y_t \mathbb{Q}_Y(dy)} \\ &\quad = \mathbb{E}^{\mathbb{Q}}\left(\mathbf{D}_t \ln\xi(Y) - \Lambda_t^{-1}\int_t^T \Lambda_s \mathbf{\Omega}_{X^*_s}(\mathbf{D}_s \ln\xi(Y))\, ds \Big| \mathcal{F}_t\right), \qquad t < T, \end{aligned}$$

where $\Lambda$ is an $\mathcal{F}$-adapted process valued in the space of $d \times d$ invertible matrices and solves the equation

$$\Lambda_t + \int_0^t \Lambda_s \boldsymbol{\Omega}_{X_s^*} \, ds = I_d.$$

Proof.   From the proof of Proposition 3, we know that for a probability $\mathbb{Q} \in \mathcal{R}_Y(\mathbb{P})$ we have

$$\mathbb{Q} = \xi(Y)\mathbb{P}.$$

Now, from the Clark–Ocone–Bismut formula,

$$\xi(Y) = 1 + \int_0^T (\Theta_s, d\tilde{X}_s),$$

where $\Theta$ is given by

$$\Theta_s = \mathbb{E}\bigg(\mathbf{D}_s \xi(Y) - \Lambda_s^{-1} \int_s^T \Lambda_u \boldsymbol{\Omega}_{X_u^*}(\mathbf{D}_u \xi(Y)) \, du | \mathcal{F}_s\bigg).$$

Girsanov's theorem gives hence

$$(X_t^*)^{-1} = \frac{\int_{\mathcal{P}} \eta_t^y \alpha_t^y \mathbb{Q}_Y(dy)}{\int_{\mathcal{P}} \eta_t^y \mathbb{Q}_Y(dy)}$$

$$= \frac{\mathbb{E}(\mathbf{D}_t \xi(Y) - \Lambda_t^{-1} \int_t^T \Lambda_s \boldsymbol{\Omega}_{X_s^*}(\mathbf{D}_s \xi(Y)) \, ds | \mathcal{F}_t)}{\mathbb{E}(\xi(Y)|\mathcal{F}_t)}$$

$$= \mathbb{E}^{\mathbb{Q}}\bigg(\mathbf{D}_t \ln \xi(Y) - \Lambda_t^{-1} \int_t^T \Lambda_s \boldsymbol{\Omega}_{X_s^*}(\mathbf{D}_s \ln \xi(Y)) \, ds | \mathcal{F}_t\bigg).$$

The last equality stems from the Bayes formula.   $\square$

By comparing Propositions 3 and 4, we deduce hence, thanks to Bayes formula, the following very general integration by parts formula, which also characterizes our process $\alpha$.

Corollary 5.   *Under the assumptions of Propositions 3 and 4, we have*

$$(3.4) \quad \begin{aligned} \int_{\mathcal{P}} \alpha_t^y \xi(y)\mathbb{P}(Y \in dy|\mathcal{F}_t) \\ = X_t^* \mathbb{E}\bigg(\mathbf{D}_t \xi(Y) - \Lambda_t^{-1} \int_t^T \Lambda_s \boldsymbol{\Omega}_{X_s^*}(\mathbf{D}_s \xi(Y)) \, ds | \mathcal{F}_t\bigg), \qquad t < T. \end{aligned}$$

Remark 6.   Let us mention here an interesting point. If we use, formally, the formula (3.4) with $\xi = \delta_y$, $y \in \mathcal{P}$, then we obtain

$$(X_t^*)^{-1} \alpha_t^y \mathbb{P}(Y \in dy|\mathcal{F}_t) = \mathbb{E}\bigg(\mathbf{D}_t \delta_Y - \Lambda_t^{-1} \int_t^T \Lambda_s \boldsymbol{\Omega}_{X_s^*}(\mathbf{D}_s \delta_Y) \, ds | \mathcal{F}_t\bigg)(dy).$$

In the flat case ($\boldsymbol{\Omega} = 0$), this formula can be found in [17], Proposition A.1, where a Malliavin calculus for measure-valued random variables is developed.

We can now give a precise definition of Brownian motion conditioned by the functional $Y$ and show how it can be constructed from a stochastic differential equation on $\mathcal{O}(\mathbb{M})$ that is called a conditioned stochastic differential equation (cf. [1]).

DEFINITION 7. A process on $\mathbb{M}$ whose law belongs to $\mathcal{R}_Y(\mathbb{P})$ is called a Brownian motion with drift $V$ conditioned by $Y$.

Let us consider a probability measure $\nu$ on the Polish space $\mathcal{P}$ which is equivalent to $\mathbb{P}_Y$. Then there exists a predictable functional on $\mathcal{C}_m(\mathbb{M})$, say $F^\nu$, such that

$$F^\nu(t, (X_s)_{0 \le s \le t}) = \frac{\int_{\mathcal{P}} \eta_t^y \alpha_t^y \nu(dy)}{\int_{\mathcal{P}} \eta_t^y \nu(dy)}, \qquad t < T.$$

DEFINITION 8. On a filtered probability space $(\boldsymbol{\Omega}, (\mathcal{H}_t)_{0 \le t \le T}, (\beta_t)_{0 \le t \le T}, \tilde{\mathbb{P}})$ which satisfies the usual conditions, where $\beta$ is an $\mathcal{H}$-adapted $d$-dimensional linear Brownian motion, the stochastic differential equation on $\mathcal{O}(\mathbb{M})$,

$$
\begin{aligned}
(3.5) \quad U_t = U_0 &+ \int_0^t V^*(U_s)\,ds \\
&+ \sum_{i=1}^d \int_0^t H_i(U_s) \circ ((U_s)^{-1} F^\nu(s, (\pi U_u)_{0 \le u \le s})\,ds + \beta_s), \qquad s \le T,
\end{aligned}
$$

is called the conditioned stochastic differential equation associated with the conditioning $(T, Y, \nu)$.

We conclude this section with the following proposition, which is a consequence of the Yamada–Watanabe theorem (see [23], page 368) which asserts that the pathwise uniqueness property for a stochastic differential equation implies uniqueness in law.

PROPOSITION 9. *Assume that* (3.5) *enjoys the pathwise uniqueness property. Then* $(\pi U_t)_{0 \le t < T}$ *is the unique Brownian motion with drift $V$ conditioned by $Y$ such that*

$$\tilde{\mathbb{P}}(Y(\pi U) \in dy) = \nu(dy).$$

PROOF.   First, we note that there exists a unique probability measure $\mathbb{Q} \in \mathcal{R}_Y(\mathbb{P})$ such that $\mathbb{Q}(Y \in dy) = \nu(dy)$. This probability $\mathbb{Q}$ is given by

$$\mathbb{Q} = \int_{\mathcal{P}} \mathbb{P}(\cdot \mid Y = y)\nu(dy).$$

Now, actually, by means of Proposition 3 we have constructed a weak solution of (3.5) on the stochastic basis $(\mathcal{C}_m(\mathbb{M}), (\mathcal{F}_t)_{0 \le t \le T}, (W_t)_{0 \le t \le T}, \mathbb{Q})$, where

$$W_t = \tilde{X}_t - \int_0^t (X_s^*)^{-1} \frac{\int_{\mathcal{P}} \eta_s^y \alpha_s^y \mathbb{Q}_Y(dy)}{\int_{\mathcal{P}} \eta_s^y \mathbb{Q}_Y(dy)} \, ds.$$

Since, thanks to the Yamada–Watanabe theorem, the pathwise uniqueness property implies uniqueness in law, we conclude that the law of $(\pi U_t)_{0 \le t \le T}$ is $\mathbb{Q}$, which exactly means that $(\pi U_t)_{0 \le t < T}$ is the unique Brownian motion conditioned by $Y$ such that

$$\tilde{\mathbb{P}}(Y(\pi U) \in dy) = \nu(dy). \qquad \square$$

## 4. Initial enlargement of Itô's filtration in the nonflat Wiener space.
In this short section, we study, under a probability measure $\mathbb{Q} \in \mathcal{R}_Y(\mathbb{P})$, the semimartingale decomposition of the coordinate process $(X_t)_{0 \le t < T}$ in the initially enlarged filtration $\mathcal{F}^Y$, where $\mathcal{F}_t^Y$ is the $\mathbb{P}$-completion of $\bigcap_{\varepsilon > 0}(\mathcal{F}_{t+\varepsilon} \vee \sigma(Y))$. This decomposition generalizes the celebrated Jacod theorem (see [1, 18, 19, 20]) in our nonflat setting.

PROPOSITION 10.    Let $\mathbb{Q} \in \mathcal{R}_Y(\mathbb{P})$. Then under $\mathbb{Q}$, the coordinate process $(X_t)_{0 \le t < T}$ is a semimartingale in the filtration $\mathcal{F}^Y$. Moreover, the process

$$\tilde{X}_t - \int_0^t (X_s^*)^{-1} \alpha_s^Y \, ds, \qquad t < T,$$

is a Brownian motion in $\mathcal{F}^Y$ under each $\mathbb{Q} \in \mathcal{R}_Y(\mathbb{P})$.

PROOF.   For almost every $y \in \mathcal{P}$, let us consider the disintegrated probability measure $\mathbb{P}^y = \mathbb{P}(\cdot \mid Y = y)$. By the very definition of $(\eta_t^y, 0 \le t < T, y \in \mathcal{P})$, the following absolute continuity relationship holds for almost every $y \in \mathcal{P}$:

$$\mathbb{P}_{/\mathcal{F}_t}^y = \eta_t^y \mathbb{Q}_{/\mathcal{F}_t}^y, \qquad t < T.$$

Thus, as a consequence of Girsanov's theorem, the process

$$\tilde{X}_t - \int_0^t (X_s^*)^{-1} \alpha_s^y \, ds, \qquad t < T,$$

is a standard Brownian motion under the probability $\mathbb{P}^y$, which implies that

$$\tilde{X}_t - \int_0^t (X_s^*)^{-1} \alpha_s^Y \, ds, \qquad t < T,$$

is a standard Brownian motion under $\mathbb{P}$ in the enlarged filtration $\mathcal{F}^Y$ (it suffices to apply Lévy's characterization of Brownian motion). Now, we note that if $\mathbb{Q} \in \mathcal{R}_Y(\mathbb{P})$, then $\mathbb{Q}$ and $\mathbb{P}$ coincide on the events which are $\mathbb{P}$ independent of $Y$. It implies that

$$\tilde{X}_t - \int_0^t (X_s^*)^{-1} \alpha_s^Y \, ds, \qquad t < T,$$

is also a Brownian motion under the probability $\mathbb{Q}$, because under $\mathbb{P}$ this process is independent of $Y$, which means that its law under $\mathbb{P}$ and $\mathbb{Q}$ is the same. $\square$

REMARK 11. It would be interesting to have conditions on $Y$ which ensure that the process $(X_t)_{0 \leq t \leq T}$ (i.e., considered up to time $T$) is a semimartingale in the filtration $\mathcal{F}^Y$. For this, we need to show that $\mathbb{Q}$-a.s.,

$$\int_0^T \|\alpha_s^Y\| \, ds < +\infty.$$

This requires an estimate for $\alpha$ which seems to be hard to obtain in all generality. For instance, it is a direct consequence from [5], page 86, that if $Y = X_T$, then the semimartingale property holds in the enlarged filtration up to time $T$.

What is really interesting in the previous proposition is that the process

$$\tilde{X}_t - \int_0^t (X_s^*)^{-1} \alpha_s^Y \, ds, \qquad t < T,$$

is a Brownian motion in $\mathcal{F}^Y$ under *any* $\mathbb{Q} \in \mathcal{R}_Y(\mathbb{P})$. As shown in the following proposition, this property characterizes $\mathcal{R}_Y(\mathbb{P})$.

PROPOSITION 12. *Let $\mathbb{Q}$ be a probability measure on $\mathcal{F}_T$ which is equivalent to $\mathbb{P}$. If the process*

$$M_t = \tilde{X}_t - \int_0^t (X_s^*)^{-1} \alpha_s^Y \, ds, \qquad t < T,$$

*is a standard Brownian motion under $\mathbb{Q}$ in the filtration $\mathcal{F}^Y$, then $\mathbb{Q} \in \mathcal{R}_Y(\mathbb{P})$.*

PROOF. For $\mathbb{P}_Y$-a.e. $y \in \mathcal{P}$, we denote $\mathbb{Q}^y$ to be the conditional probability $\mathbb{Q}(\cdot | Y = y)$. From our assumption, the process $M$ is, under $\mathbb{Q}^y$, a standard Brownian motion. Hence, by Girsanov's theorem,

$$d\mathbb{Q}^y_{/\mathcal{F}_t} = \eta_t^y \, d\mathbb{P}_{/\mathcal{F}_t}, \qquad t < T.$$

Since we also have

$$d\mathbb{P}^y_{/\mathcal{F}_t} = \eta^y_t d\mathbb{P}_{/\mathcal{F}_t}, \qquad t < T,$$

where $\mathbb{P}^y$ is the conditional probability $\mathbb{P}(\cdot|Y = y)$, we immediately deduce

$$\mathbb{Q}^y = \mathbb{P}^y$$

and hence $\mathbb{Q} \in \mathcal{R}_Y(\mathbb{P})$.   $\square$

## 5. Examples.

5.1. *The case $Y = X_T$ and the infinite pinning class.* The first example we can think of is the example $Y = X_T$, $T > 0$. It means that we condition the law at $T$ of the Brownian motion with drift $V$. In what follows, for a smooth function $f(x, y)$, $x, y \in \mathbb{M}$, the vector field $\nabla f$ denotes the gradient computed with respect to the first variable and $Q_t$ denotes the semigroup associated with the transition function $q_t$ (cf. Section 2 for the notations).

In this special case, it is easily seen that the assumption (3.1) is satisfied with

$$\eta^y_t = \frac{q_{T-t}(X_t, y)}{q_T(m, y)}, \qquad t < T, y \in \mathbb{M}.$$

Moreover, a direct computation based on Itô's formula shows that

$$\alpha^y_t = \nabla \ln q_{T-t}(X_t, y), \qquad t < T, y \in \mathbb{M}.$$

The formula for the compensator of $\tilde{X}$ under $\mathbb{Q} \in \mathcal{R}_{X_T}(\mathbb{P})$ is hence

$$(X^*_t)^{-1} \frac{\int_{\mathbb{M}} \eta^y_t \alpha^y_t \mathbb{Q}_Y(dy)}{\int_{\mathbb{M}} \eta^y_t \mathbb{Q}_Y(dy)} = (X^*_t)^{-1} \nabla \ln Q_{T-t} \xi(X_t), \qquad t < T,$$

where $\xi := d\mathbb{Q}_Y/d\mathbb{P}_Y$. Notice that from this, we deduce, by rewriting Proposition 4, that for $\mathbb{Q} \in \mathcal{R}_{X_T}(\mathbb{P})$,

$$(5.1) \quad \nabla \ln Q_{T-t} \xi(X_t) = \mathbb{E}^{\mathbb{Q}}(X^*_t \Lambda^{-1}_t \Lambda_T (X^*_T)^{-1} (\nabla \ln \xi)(X_T)|\mathcal{F}_t), \qquad t < T,$$

and that formula (3.4) reads, in this case,

$$\nabla Q_{T-t} \xi(X_t) = \mathbb{E}(X^*_t \Lambda^{-1}_t \Lambda_T (X^*_T)^{-1} (\nabla \xi)(X_T)|\mathcal{F}_t), \qquad t < T,$$

which is a fairly well-known formula (see e.g., [8] and [27]).

Let us mention a corollary to this. In [4], it was shown that if $(\tilde{Z}_t)_{0 \le t \le T}$ denotes the solution of the stochastic differential equation

$$\tilde{Z}_t = \int_0^t b(\tilde{Z}_s) \, ds + B_t, \qquad 0 \le t \le T,$$

where $B$ is one-dimensional linear Brownian motion and $b$ is a smooth function, then the process

$$\left( \exp\left( \int_0^t b'(\tilde{Z}_s)\,ds \right) \frac{\partial}{\partial x} \ln \tilde{q}_{T-t}(\tilde{Z}_t, \tilde{Z}_T) \right)_{0 \le t < T}$$

is a martingale in the natural filtration of $\tilde{Z}$ initially enlarged by $\tilde{Z}_T$; $\tilde{q}_{T-t}$ denotes here the transition function of $\tilde{Z}$. In our nonflat setting, we have the following analogue of this result.

COROLLARY 13. *The process*

$$(\Lambda_t(X_t^*)^{-1} \nabla \ln q_{T-t}(X_t, X_T))_{0 \le t < T}$$

*is a martingale in the enlarged filtration $\mathcal{F}^{X_T}$ under each probability $\mathbb{Q} \in \mathcal{R}_{X_T}(\mathbb{P})$.*

PROOF. It is enough to show that under $\mathbb{P}^y$ the process

$$(\Lambda_t(X_t^*)^{-1} \nabla \ln q_{T-t}(X_t, y))_{0 \le t < T}$$

is a martingale, where $\mathbb{P}^y$ is the disintegrated probability $\mathbb{P}(\cdot | X_T = y)$, because for $\mathbb{Q} \in \mathcal{R}_{X_T}(\mathbb{P})$,

$$\mathbb{Q} = \int_{\mathbb{M}} \mathbb{P}^y \mathbb{Q}_{X_T}(dy).$$

Let $\xi_n$ be a sequence of smooth positive and normalized functions such that

$$\xi_n(x) q_T(m, x)\,dx \xrightarrow[n \to +\infty]{weakly} \delta_y,$$

where $\delta_y$ denotes the Dirac measure at $y$. As a direct consequence of (5.1) we deduce that

$$\left( \Lambda_t(X_t^*)^{-1} \frac{\int_{\mathbb{M}} \nabla q_{T-t}(X_t, y) \xi_n(y)\,dy}{\int_{\mathbb{M}} q_{T-t}(X_t, y) \xi_n(y)\,dy} \right)_{0 \le t < T}$$

is a martingale under $\mathbb{Q}^n = \xi_n(Y_T)\,\mathbb{P}$. Since

$$\mathbb{Q}^n \xrightarrow[n \to +\infty]{weakly} \mathbb{P}^y,$$

we conclude that

$$(\Lambda_t(X_t^*)^{-1} \nabla \ln q_{T-t}(X_t, y))_{0 \le t < T}$$

is a martingale under $\mathbb{P}^y$. $\square$

Notice that the previous martingale can be used to recover the generalization from [10] (see also [8]) of the celebrated Bismut formula (see [5], formula 2.80, page 283), which express $\nabla \ln p_{T-t}(X_t, y)$ as the expectation of a stochastic integral, where $p_t$ denotes the heat transition function. Indeed, since the process

$$(\Lambda_t(X_t^*)^{-1} \nabla \ln q_{T-t}(X_t, X_T))_{0 \leq t < T}$$

is martingale in the enlarged filtration, we deduce

$$\nabla \ln q_T(m, X_T) = \mathbb{E}^{\mathbb{Q}}(\Lambda_t(X_t^*)^{-1} \nabla \ln q_{T-t}(X_t, X_T)|X_T).$$

By taking the expectation under $\mathbb{Q}$ on the both sides of the equality, we deduce

$$\nabla \ln Q_T \xi(m) = \mathbb{E}^{\mathbb{Q}}(\Lambda_t(X_t^*)^{-1} \nabla \ln q_{T-t}(X_t, X_T)).$$

Now, it suffices to note that

$$\frac{1}{T}\mathbb{E}^{\mathbb{Q}}\left(\int_0^T \Lambda_s \, d\tilde{X}_s\right) = \frac{1}{T}\mathbb{E}^{\mathbb{Q}}\left(\int_0^T \Lambda_s(X_s^*)^{-1} \nabla \ln q_{T-s}(X_s, X_T) \, ds\right)$$

to obtain

$$\nabla \ln Q_T \xi(m) = \frac{1}{T}\mathbb{E}^{\mathbb{Q}}\left(\int_0^T \Lambda_s \, d\tilde{X}_s\right),$$

which is Theorem 3.1 of [10] (see also [8] and [27]).

Now, we are interested in the infinite pinning class of the Brownian motion with drift $V$. We define this set of probabilities on $\mathcal{C}_m(\mathbb{M})$ by

$$\mathcal{R}^\infty = \bigcap_{T>0} \mathcal{R}_{X_T}(\mathbb{P}).$$

Notice that a probability which belongs to $\mathcal{R}^\infty$ is not necessarily equivalent or even absolutely continuous with respect to $\mathbb{P}$. A process $(\tilde{Z}_t)_{t \geq 0}$ on $\mathbb{M}$ whose law belongs to $\mathcal{R}^\infty$ is a nonhomogeneous diffusion which has the same bridges as $(Z_t)_{t \geq 0}$ (this process is defined in Section 2). Precisely for any $T > 0$ and $n \in \mathbb{M}$,

$$\mathbb{P}(\cdot|\tilde{Z}_T = n) = \mathbb{P}(\cdot|Z_T = n).$$

In the case $\mathbb{M} = \mathbb{R}$, these processes were characterized in [13] (see also [1]).

Let us recall that a Borel function $\phi: \mathbb{R}_+ \times \mathbb{M} \to \mathbb{R}$ is called a $\mathbb{P}$-space–time harmonics if the process $(\phi(t, X_t))_{t \geq 0}$ is $\mathbb{P}$-martingale. The following easy proposition shows that there is a bijection between $\mathcal{R}^\infty$ and the set of continuous $\mathbb{P}$-space–time positive and normalized harmonics.

PROPOSITION 14. *Let $\phi$ be continuous $\mathbb{P}$-space–time positive harmonics. Then there is a unique $\mathbb{Q} \in \mathcal{R}^\infty$ such that for any $T > 0$,*

$$\mathbb{Q}_{/\mathcal{F}_T} = \frac{\phi(T, X_T)}{\phi(0, m)} \mathbb{P}_{/\mathcal{F}_T}.$$

*Conversely, if $\mathbb{Q} \in \mathcal{R}^\infty$, then there is a unique continuous $\mathbb{P}$-space–time positive harmonics $\phi$ such that for any $T > 0$,*

$$\mathbb{Q}_{/\mathcal{F}_T} = \phi(T, X_T)\mathbb{P}_{/\mathcal{F}_T}.$$

We deduce the following corollary.

COROLLARY 15. *Assume that $\mathbb{M}$ is compact and that there exists a constant $k > 0$ such that for any vector field $Z$ on $\mathbb{M}$ such that $\|Z\| \leq 1$,*

$$(5.2) \qquad \tfrac{1}{2}\mathbf{Ric}(Z, Z) - (\nabla_Z V, Z) \geq k\|Z\|^2.$$

*Then $\mathcal{R}^\infty = \{\mathbb{P}\}$.*

PROOF. Let $\mathbb{Q} \in \mathcal{R}^\infty$. There exists a unique continuous $\mathbb{P}$-space–time positive harmonics $\varphi$ such that for any $T > 0$,

$$\mathbb{Q}_{/\mathcal{F}_T} = \varphi(T, X_T)\mathbb{P}_{/\mathcal{F}_T}.$$

Being a $\mathbb{P}$-space–time harmonics, $\varphi$ is a weak solution of the backward heat equation

$$\frac{\partial \varphi}{\partial t} + \frac{1}{2}\Delta\varphi + V\varphi = 0.$$

Moreover it is bounded on $[0, T] \times \mathbb{M}$ by the compactness of $\mathbb{M}$. Thus, actually $\varphi$ is smooth and is a strong solution.

For $T > 0$, the process

$$(\Lambda_t(X_t^*)^{-1}\nabla\ln q_{T-t}(X_t, X_T))_{0 \leq t < T}$$

is a martingale in the enlarged filtration $\mathcal{F}^{X_T}$ under the probability $\mathbb{Q}$. This implies that for $t \geq 0$,

$$\Lambda_s(X_s^*)^{-1}\nabla\ln\varphi(s, X_s) = \mathbb{E}^{\mathbb{Q}}(\Phi_t(X_t^*)^{-1}\nabla\ln q_1(X_t, X_{t+1})|\mathcal{F}_s).$$

Now, the assumption (5.2) implies, thanks to Gronwall's lemma, that

$$\|\Lambda_t\| \leq e^{-kt}.$$

Hence, since $\nabla\ln q_1$ is bounded because $\mathbb{M}$ is compact, there exists a constant $K$ such that

$$\|\Lambda_s\|\|\nabla\ln\varphi(s, X_s)\| \leq Ke^{-kt}$$

by letting $t \to +\infty$. We deduce $\|\nabla\ln\varphi(s, X_s)\| = 0$, which implies $\varphi = 1$ and $\mathbb{Q} = \mathbb{P}$. □

5.2. *Conditioning the first hitting time of a small geodesic sphere.* In this section, we assume that $\mathbb{M}$ is compact and that $V = 0$. In this case, $\mathbb{P}$ is simply the Wiener measure. We consider here the small geodesic sphere with radius $r > 0$ and centered at $m$. By small, we mean that $r$ is lower than the injectivity radius at $m$. We denote this sphere by $\mathbb{S}_r(m)$. Although the random variable

$$T_r = \inf\{t > 0, X_t \in \mathbb{S}_r(m)\}$$

does not satisfy the assumption (3.1), there exists a process

$$(\eta_t^\tau, t < T_r, \tau > 0)$$

such that

$$\mathbb{P}(T_r \in d\tau | \mathcal{F}_t, t < T_r) = \eta_t^\tau \mathbb{P}(T_r \in d\tau).$$

Namely we have

$$\eta_t^\tau = \Psi^\tau(t, X_t),$$

where $\Psi^\tau$, $\tau \in \mathbb{R}_+^*$, is the solution on the open ball centered at $m$ with radius $r$ of the terminal value problem

$$\frac{1}{2}\Delta\Psi^\tau + \frac{\partial \Psi^\tau}{\partial t} = 0,$$

$$\Psi^\tau(0, m) = 1,$$

$$\Psi^\tau_{/\mathbb{S}_r(m)} = \delta_\tau,$$

where $\delta_\tau$ is the Dirac distribution at $\tau$. This allows us to use the results of Sections 2 and 3 (cf. Remark 1) up to the stopping time $T_r$; in particular, a process

$$(\alpha_t^\tau, t < T_r, \tau > 0)$$

can be defined in the same way as in Lemma 2. Precisely, we have

$$\alpha_t^\tau = \nabla \ln \Psi^\tau(t, X_t).$$

The formula for the compensator of $\tilde{X}$ under $\mathbb{Q} \in \mathcal{R}_{T_r}(\mathbb{P})$ is hence

$$(X_t^*)^{-1} \frac{\int_{\mathbb{R}_+^*} \eta_t^\tau \alpha_t^\tau \mathbb{Q}_Y(d\tau)}{\int_{\mathbb{R}_+^*} \eta_t^\tau \mathbb{Q}_Y(d\tau)} = (X_t^*)^{-1}\nabla \ln \varphi(t, X_t), \qquad t < T_r,$$

where $\varphi$ is the solution of the terminal value problem

$$\frac{1}{2}\Delta\varphi + \frac{\partial \varphi}{\partial t} = 0,$$

$$\varphi(0, m) = 1,$$

$$\varphi_{/\mathbb{S}_r(m)} = \frac{d\mathbb{Q}_{T_r}}{d\mathbb{P}_{T_r}}.$$

By space–time duality for Markov processes, it is interesting to note that we also have a martingale associated with the pinning class $\mathcal{R}_{T_r}(\mathbb{P})$. Actually, this martingale is constructed exactly as in Corollary 13, but by reasoning up to time $T_r$.

PROPOSITION 16. *The process*

$$(\Phi_t(X_t^*)^{-1}\nabla \ln \Psi^{T_r}(t, X_t))_{0 \leq t < T_r}$$

*is a martingale in the enlarged filtration $\mathcal{F}^{T_r}$ under each probability $\mathbb{Q} \in \mathcal{R}_{T_r}(\mathbb{P})$, where $\Phi$ is an adapted process valued in the space of $d \times d$ invertible matrices and solves the equation*

$$\Phi_t + \frac{1}{2}\int_0^t \Phi_s \overline{\mathbf{Ric}}_{X_s^*}\, ds = I_d, \qquad t < T_r.$$

To conclude this paper, we show now that the pinning classes $\mathcal{R}_{T_r}^m$, $m \in \mathbb{M}$, characterize the isotropy of the manifold (the compactness of $\mathbb{M}$ is necessary).

PROPOSITION 17. *Let us assume that $\mathbb{M}$ is simply connected and that for all $m \in \mathbb{M}$ there exists a small $\varepsilon > 0$ such that*

$$\bigcap_{0 < r < \varepsilon} \mathcal{R}_{T_r}^m(\mathbb{P}) \supsetneq \{\mathbb{P}_m\}.$$

*Then, $\mathbb{M}$ is isometric to a compact rank 1 symmetric space.*

PROOF. The assumption of the nontriviality of the pinning classes

$$\bigcap_{0 < r < \varepsilon} \mathcal{R}_{T_r}^m(\mathbb{P})$$

implies that for all $m \in \mathbb{M}$, there exists a harmonic radial function $n \to \varphi(\delta(m, n))$ defined on a neighborhood of $m$ ($\delta$ denotes the Riemannian distance on $\mathbb{M}$). From the Lichnerowicz–Szabo theorem (see [26]; we also refer to [2] for a probabilistic understanding and proof of this theorem), this implies that $\mathbb{M}$ is isometric to a compact rank 1 symmetric space. □

**6. Opening.** In the linear case, it was shown in [4] that the one-parameter group of the translations on the path space

$$T_\alpha : \mathcal{C}_0(\mathbb{R}^d) \to \mathcal{C}_0(\mathbb{R}^d),$$

$$\omega \to (\omega_t + \alpha t)_{t \geq 0}$$

acts naturally on the pinning class $\mathcal{R}_{X_T}(\mathbb{P})$. This action is closely related to the quasi-invariance of the Wiener measure by the translations. So, it seems to us that it would be interesting to study the one-parameter groups of adapted transformations which act naturally on the nonflat pinning class.

**Acknowledgments.** This work was performed while I was visiting the University of Barcelona from September 6, 2002 until December 6, 2002. I cordially thank David Nualart and Marta Sanz for their hospitality. I also thank anonymous referees for judicious comments.

## REFERENCES

[1] BAUDOIN, F. (2002). Conditioned stochastic differential equations: Theory, examples and applications to finance. *Stochastic Process. Appl.* **100** 109–145. MR1919610

[2] BAUDOIN, F. (2002). Skew-product decompositions of Brownian motions on manifolds: A probabilistic aspect of Lichnerowicz–Szabo theorem. *Bull. Sci. Math.* **126** 481–491. MR1931625

[3] BAUDOIN, F. (2003). *Modelling Anticipations on a Financial Market. Lectures Notes in Math.* **1814**. Springer, Berlin. MR2021790

[4] BAUDOIN, F. and THIEULLEN, M. (2003). Pinning class of the Wiener measure by a functional: Related martingales and invariance properties. *Probab. Theory Related Fields* **127** 1–36. MR2006229

[5] BISMUT, J. M. (1984). *Large Deviations and the Malliavin Calculus.* Birkhäuser, Boston. MR755001

[6] BISMUT, J. M. (1984). The Atiyah–Singer theorems: A probabilistic approach, I. The index theorem. *J. Funct. Anal.* **57** 56–99. MR744920

[7] BISMUT, J. M. (1984). The Atiyah–Singer theorems: A probabilistic approach, II. The Lefschetz fixed point formulas. *J. Funct. Anal.* **57** 329–348. MR756173

[8] DRIVER, B. and THALMAIER, A. (2001). Heat equation derivative formulas for vector bundles. *J. Funct. Anal.* **1** 42–108. MR1837533

[9] ELWORTHY, K. D. (1982). *Stochastic Differential Equations on Manifolds.* Cambridge Univ. Press. MR675100

[10] ELWORTHY, K. D. and LI, X. M. (1994). Formulae for the derivatives of heat semigroups. *J. Funct. Anal.* **125** 252–286. MR1297021

[11] EMERY, M. (1989). *Stochastic Calculus in Manifolds.* Springer, Berlin. MR1030543

[12] EMERY, M. (2000). Martingales continues dans les variétés différentiables. *Ecole d'Eté de Saint-Flour XXVIII. Lecture Notes in Math.* **1738** 1–84. MR1775639

[13] FITZSIMMONS, P. (1998). Markov processes with identical bridges. *Electron. J. Probab.* **3**. MR1641066

[14] FITZSIMMONS, P., PITMAN, J. W. and YOR, M. (1993). Markovian bridges, construction, Palm interpretation and splicing. In *Seminar on Stochastic Processes* (E. Çinlar, K. L. Chung and M. J. Sharpe, eds.) 101–134. Birkhäuser, Boston. MR1278079

[15] HSU, E. P. (2002). *Stochastic Analysis on Manifolds.* Springer, New York. MR1882015

[16] IKEDA, N. and WATANABE, S. (1989). *Stochastic Differential Equations and Diffusion Processes*, 2nd ed. North-Holland, Amsterdam. MR1011252

[17] IMKELLER, P., PONTIER, M. and WEISZ, F. (2001). Free lunch and arbitrage possibilities in a financial market with an insider. *Stochastic Process. Appl.* **92** 103–130. MR1815181

[18] JACOD, J. (1985). Grossissement initial, hypothèse $(H')$ et théorème de Girsanov. *Grossissements de filtrations*: *Exemples et applications. Lecture Notes in Math.* **1118** 15–35. Springer, Berlin. MR883648

[19] JEULIN, T. (1980). *Semi-Martingales et Grossissement d'une Filtration. Lecture Notes in Math.* **833**. Springer, Berlin. MR604176

[20] Jeulin, T. and Yor, M., eds. (1985). *Grossissements de Filtrations*: *Exemples et Applications. Lecture Notes in Math.* **1118**. Springer, Berlin. MR884713

[21] Malliavin, P. (1997). *Stochastic Analysis.* Springer, Berlin. MR1450093

[22] Nualart, D. (1995). *The Malliavin Calculus and Related Topics.* Springer, Berlin. MR1344217

[23] Revuz, D. and Yor, M. (1999). *Continuous Martingales and Brownian Motion*, 3rd ed. Springer, Berlin. MR1725357

[24] Stroock, D. W. (2000). *An Introduction to the Analysis of Paths on Riemannian Manifolds.* Amer. Math. Soc., Providence, RI. MR1715265

[25] Stroock, D. W. and Varadhan, S. R. S. (1979). *Multidimensional Diffusion Processes.* Springer, Berlin. MR532498

[26] Szabo, S. (1990). The Lichnerowicz conjecture on harmonic manifolds. *J. Differential Geom.* **31** 1–28. MR1030663

[27] Thalmaier, A. (1997). On the differentiation of heat semigroups and Poisson integrals. *Stochastics Stochastics Rep.* **61** 297–321. MR1488139

[28] Thieullen, M. and Zambrini, J. C. (1997). Symmetries in the stochastic calculus of variations. *Probab. Theory Related Fields* **107** 401–427. MR1440139

[29] Yor, M. (1997). *Some Aspects of Brownian Motion, Part II.* Birkhäuser, Basel. MR1442263

Laboratoire de Probabilités
et Statistiques
Université Paul Sabatier
118 route de Narbonne 31062
Toulouse Cedex 4
France
e-mail: symplectik@aol.com
e-mail: fbaudoin@cict.fr